\documentclass[12pt]{article}

\newcommand{\mk}{\ensuremath{m_k }} 
\newcommand{\ml}{\ensuremath{m_{\ell} }} 

\newcommand{\Kdd}{\ensuremath{{K_{d,d}}}} 
\newcommand{\DK}{\ensuremath{DK_{N,d}\,}} 

\def\BPi{\displaystyle\prod}  

\def\PFG{\ensuremath{Z^{\rm{match}}_{\lambda}(G)}} 

\def\PFG{\ensuremath{Z^{\rm{match}}_{\lambda}(G)}} 

\def\qed{\hfill$\square$\bigskip\medskip}
\def\sqr#1#2{{\vcenter{\vbox{\hrule height.#2pt
\hbox{\vrule width.#2pt height #1pt \kern#1pt \vrule width.#2pt}
\hrule height.#2pt}}}}

\newtheorem{newthm}{}[section]
\newcommand{\comment}[1]{}

\newtheorem{theorem}[newthm]{Theorem}
\newtheorem{lemma}[newthm]{Lemma}

\newtheorem{conjecture}[newthm]{Conjecture}

\newtheorem{thm}{Theorem}[section]

\newtheorem{remark}[thm]{Remark}

\newcommand{\cE}[0]{{\cal E}}

\newcommand{\cI}[0]{{\cal I}}

\newcommand{\cO}[0]{{\cal O}}

\newcommand{\ga}[0]{\alpha}

\newcommand{\gl}[0]{\lambda}

\usepackage{amsmath,amssymb}
\usepackage{mathrsfs}

\title{Matchings and Independent Sets of a Fixed Size in Regular Graphs}
\author{Teena Carroll\thanks{Mathematics Department, St. Norbert College, De Pere WI.}~~~~David Galvin\thanks{Department of Mathematics, University of
Notre Dame, South Bend IN}~~~~Prasad Tetali\thanks{School of
Mathematics \& School of Computer Science, Georgia Institute of
Technology, Atlanta GA. Research supported in part by NSF
grant DMS-0701043.}}
\date{November 2008}

\begin{document}
\renewcommand{\thefootnote}{\fnsymbol{footnote}}
\bibliographystyle{abbrv}
\maketitle

\begin{abstract}
We use an entropy based method to study two 
graph maximization problems.
We upper bound the number of matchings of fixed size $\ell$ in a $d$-regular graph on $N$ vertices. For $\frac{2\ell}{N}$ bounded away from $0$ and $1$, the logarithm of the bound we obtain agrees in its leading term with the logarithm of the number of matchings of size $\ell$ in the graph consisting of $\frac{N}{2d}$ disjoint copies of $\Kdd$. This provides asymptotic evidence for a conjecture of S. Friedland {\it et al.}.
We also obtain an analogous
result for independent sets of a fixed size in regular graphs, giving asymptotic evidence for a conjecture of J. Kahn.  Our bounds on the number of matchings and independent sets of a fixed size are derived from bounds on the partition function (or generating polynomial) for matchings and independent sets.
\end{abstract}

\section{Introduction}

Given a $d$-regular graph $G$ on $N$ vertices and a particular type
of subgraph, a natural class of problems arises:  ``How many
subgraphs of this type can $G$ contain?''  In this paper we give
upper bounds on the number of partial matchings of a fixed
fractional size, and on the number of independent sets of a fixed
size, in a general $d$-regular graph, and we show that our bounds
are asymptotically matched at the logarithmic level by the graph
consisting of $\frac{N}{2d}$ disjoint copies of $\Kdd$. (See
\cite{Bollo1998} and \cite{Diest2005} for graph theory basics.)

Let $G$ be a bipartite graph on $N$ vertices with partition classes
$A$ and $B$ and with $|A|=|B|$. Suppose that the degree sequence of $A$ is given
by $\{r_i\}_{i=1}^{|A|}$. A result of Br\'{e}gman concerning the permanent of $0$-$1$ matrices
\cite{Bregm1973} (see also \cite{Alon2000}) gives a bound on the
number of perfect matchings in $G$:
\begin{theorem}{
\label{Bregman} (Br\'{e}gman) Let $\mathcal{M}_{\rm{perfect}}(G)$ be
the set of perfect matchings in $G$.  Then
\[ \left|\mathcal{M}_{\rm{perfect}}(G)\right| \leq \BPi_{i=1}^{|A|}(r_i!)^{\frac{1}{r_i}}. \]
}
\end{theorem}

When $r_i=d$ for all $i$ and $|A|$ is divisible by $d$, equality in
the above theorem is achieved by the graph consisting of
$\frac{N}{2d}$ disjoint copies of  the complete bipartite graph
\Kdd, so we know that among $d$-regular bipartite graphs on $N$
vertices, with $2d|N$, this graph contains the greatest number of
perfect matchings. 
(Wanless \cite{Wanless2003} has considered the case when $2d$ is not a multiple of $N$, obtaining lower bounds on $\left|\mathcal{M}_{\rm{perfect}}(G)\right|$ and some structural results on the maximizing graphs in this case.)

\comment{
Friedland \textit{et al.} \cite{Friedland2005} propose an extension
of this observation, which they call the ``upper matching
conjecture.'' They conjecture that among all $d$-regular bipartite
graphs on $N$ vertices, with $2d|N$, none has more matchings of size
$\ell$ for $0 \leq \ell \leq \frac{N}{2}$ than the graph consisting
of $\frac{N}{2d}$ copies of \Kdd. }

Friedland \textit{et al.} \cite{Friedland2005} propose an extension
of this observation, which they call the Upper Matching
Conjecture. Write $\ml(G)$ for the number of matchings in $G$ of size $\ell$, and write $\DK$ for the graph consisting of
$\frac{N}{2d}$ disjoint copies of \Kdd.
\begin{conjecture} \label{conj-UMC}
For any $N$-vertex, $d$-regular graph $G$ with $2d|N$ and any $0
\leq \ell \leq N/2$,
$$
\ml(G) \leq \ml(\DK).
$$
\end{conjecture}

%
%

In this note we upper bound the logarithm of the number of $\ell$-matchings of a
regular graph and show that, at the level of the leading term, this
upper bound is achieved by the disjoint union of the appropriate
number of copies of \Kdd$\,$.
We will use the parameterization
%
$\alpha=\frac{2\ell}{N}$, and
refer interchangeably to a matching of size $\ell$ or a matching
whose size is an $\alpha$-fraction of the maximum possible matching
size.
%
In what follows, $H(x)=-x\log x-(1-x)\log(1-x)$ is the usual binary
entropy function. (All logarithms in this note are base 2.)

\begin{theorem} \label{thm-main_match}
\label{dl} {Let $G$ be a $d$-regular graph on $N$ vertices and
$\ell$ an integer satisfying $0\leq \ell \leq \frac{N}{2}$. Set
$\alpha=\frac{2\ell}{N}$. The number of matchings in $G$ of size
$\ell$ satisfies
\[    \log(\ml(G)) \leq \frac{N}{2}\big[\alpha \log d+H(\alpha)].      \]
This bound is tight up to the first order term: for fixed $\alpha
\in (0,1)$,
\[ \log(\ml(\DK))\geq \frac{N}{2}\left[\alpha\log d+2H(\alpha)+
\alpha\log\left(\frac{\alpha}{e}\right)+\Omega\left(\frac{\log d}{d}\right)\right],
\]
with the constant in the $\Omega$ term depending on $\alpha$.}
\end{theorem}

In \cite{Friedland2008} an asymptotic variant of Conjecture \ref{conj-UMC} is presented. Let $\{G_k\}$ be a sequence of $d$-regular bipartite
graphs with $|V_k|$, the number of vertices of $G_k$, growing to infinity, and fix $\alpha \in [0,1]$. Set
$$
h_{\{G_k\}}(\alpha) = \limsup (\log m_{\ell_k}(G_k))/|V_k|
$$
where the limit is over all sequences $\{\ell_k\}$ with
$2\ell_k/|V_k| \rightarrow \alpha$. The Asymptotic Upper Matching Conjecture asserts that
$$
h_{\{G_k\}}(\alpha) \leq h_{\{kK_{d,d}\}}(\alpha)
$$
where $kK_{d,d}$ is the graph consisting of $k$ disjoint copies of $K_{d,d}$. Theorem \ref{thm-main_match} shows that for each fixed $\alpha$, there is a constant $c_\alpha$ (independent of $d$) with $h_{\{G_k\}}(\alpha) \leq h_{\{kK_{d,d}\}}(\alpha) + c_\alpha$.

%

\medskip

We show similar results for the number of independent sets in
$d$-regular graphs.  A point of departure for our consideration of
independent sets is the following result of Kahn \cite{Kahn2001}.
For any graph $G$ write $\cI(G)$ for the set of independent sets in
$G$ and write $i_t(G)$ for the set of independent sets of size $t$
(i.e., with $t$ vertices).

\begin{theorem}{
\label{thm-ind_sets} (Kahn) For any $N$-vertex, $d$-regular
bipartite graph $G$,
$$
|\cI(G)| \leq |\cI(K_{d,d})|^{N/2d}.
$$}
\end{theorem}
Note that when $2d|N$, we have $|\cI(K_{d,d})|^{N/2d}=|\cI(\DK)|$.
Kahn \cite{Kahn2001} proposes the following natural conjecture.
\begin{conjecture} \label{conj-ind}
For any $N$-vertex, $d$-regular graph $G$ with $2d|N$ and any $0
\leq t \leq N/2$,
$$
i_t(G) \leq i_t(\DK).
$$
\end{conjecture}
We provide asymptotic evidence for this conjecture.
\begin{theorem} \label{thm-indsets}
For $N$-vertex, $d$-regular $G$, and $0 \leq t \leq N/2$,
\begin{equation} \label{ind_set_ub}
i_t(G) \leq \left\{
\begin{array}{ll}
2^{\frac{N}{2}\left(H\left(\frac{2t}{N}\right)+\frac{2}{d}\right)} &
\mbox{in general} \\
2^{\frac{N}{2}\left(H\left(\frac{2t}{N}\right)+\frac{1}{d} -
\frac{\log e}{2d}\left(1-\frac{2t}{N}\right)^d \right)} & \mbox{if
$G$ is bipartite} \\
2^t{\frac{N}{2} \choose t} & \mbox{if $G$ has a perfect matching.}
\end{array}
\right.
\end{equation}
On the other hand,
\begin{equation} \label{ind_set_lb}
i_t(\DK) \geq \left\{
\begin{array}{ll}
\left(1-\frac{1}{c}\right)\binom{\frac{N}{2}}{t}2^{\frac{N}{2}\left(\frac{1}{d}-\frac{c}{d}\left(1-\frac{2t}{N}\right)^d\right)}
& \mbox{for any $c > 1$} \\
& \\
2^t{\frac{N}{2} \choose t} \prod_{k=1}^{t-1}
\left(1-\frac{2kd}{N}\right) & \mbox{for $t \leq \frac{N}{2d}$.}
\end{array}
\right.
\end{equation}
\end{theorem}


If $N$, $d$ and $t$ are sequences satisfying $t=\alpha
\frac{N}{2}$ for some fixed $\alpha \in (0,1)$ and $G$ is a sequence
of $N$-vertex, $d$-regular graphs, then from
(\ref{ind_set_ub})
$$
\log i_t(G) \leq \left\{
\begin{array}{ll}
\frac{N}{2}\left[H\left(\alpha\right)+\frac{2}{d}\right] &
\mbox{in general} \\
& \\
 \frac{N}{2}\left[H\left(\alpha\right)+\frac{1}{d}\right] &
\mbox{if $G$ is bipartite,}
\end{array}
\right.
$$
whereas if $N=\omega(d \log d)$ and $d=\omega(1)$ then taking $c=2$
in the first bound of (\ref{ind_set_lb}) and using Stirling's
formula to analyze the behavior of ${N/2 \choose \alpha N/2}$, we
obtain the near matching lower bound
$$
\log i_t(\DK)  \geq
\frac{N}{2}\left[H\left(\alpha\right)+\frac{1}{d}(1+o(1))\right].
$$

If $N=o\left(d/(1-\alpha)^d\right)$ and $G$ is bipartite, then the
gap between our bounds on $i_t(G)$ and $i_t(\DK)$ is just a multiplicative factor of $O(\sqrt{N})$; indeed, in this
case (taking any $c=\omega(1)$) we obtain from the first bound of
(\ref{ind_set_lb}) that
$$
i_t(\DK) \geq
(1-o(1))\binom{\frac{N}{2}}{t}2^{\frac{N}{2}\left(H(\ga)+\frac{1}{d}\right)}.
$$

For smaller sets, whose sizes scale with $N/d$ rather than $N$,
the final bounds in (\ref{ind_set_ub}) and (\ref{ind_set_lb}) come
into play. Specifically, for any $N$, $t$ and $d$
\begin{equation} \label{ind_set_lb_spec}
i_t(\DK) \geq \left\{
\begin{array}{ll}
{\frac{N}{2} \choose t}2^{t\left(1+o(1)\right)} & \mbox{if
$t=o\left(\frac{N}{d}\right)$} \\
(1+o(1)){\frac{N}{2} \choose t}2^t & \mbox{if
$t=o\left(\sqrt{\frac{N}{d}}\right)$}
\end{array}
\right.
\end{equation}
Note that in the latter case, for $G$ with a perfect matching we
have $i_t(G) \leq (1+o(1))i_t(\DK)$. To obtain
(\ref{ind_set_lb_spec}) from (\ref{ind_set_lb}) we use
$$
\prod_{k=1}^{t-1}\left(1-\frac{2kd}{N}\right) \geq
\exp\left\{-\frac{4d}{N}\sum_{k=1}^{t-1}k\right\} \geq \exp\left\{-\frac{2dt(t-1)}{N}\right\}.
$$

\section{Counting Matchings}
Given a graph $G$ and a nonnegative real number $\lambda$, we can
form weighted matchings of $G$ by assigning each matching containing
$\ell$ edges weight $\lambda^{\ell}$. The weighted partition
function, $\PFG$, gives the total weight of matchings.  Formally,
$$
\PFG:= \sum_{m \in {\cal M}(G)}
\lambda^{|m|}=\sum_{k=0}^{\frac{N}{2}} \mk(G) \, \lambda^k.
$$

(This is often referred to as the generating function for matchings or the matching polynomial). We will prove Theorem \ref{dl} by showing a bound on the partition function, and then using
that bound to limit the number of matchings of a particular weight (size).
\begin{lemma}
\label{pfg} {For all $d$-regular graphs $G$, $\PFG \leq
(1+d\lambda)^{\frac{N}{2}}$}
\end{lemma}

This lemma is easily proven in the bipartite case; the difficulty
arises when we want to prove the same bound for general graphs.
Indeed, if $G$ is a bipartite graph with bipartition classes $A$ and
$B$, we can easily see that the right hand side above counts a
superset of weighted matchings.  Elements in this superset are sets
of edges no two of which are adjacent to the same element of $A$
(but with no restriction on incidences with $B$).

\medskip

\noindent \textbf{Proof of Lemma \ref{pfg}} To prove this lemma, we
will use the following result of Friedgut  \cite{Friedgut2004},
which describes a weighted version of the information theoretic
Shearer's Lemma.

\begin{theorem}(Friedgut)
\label{weighted_entropy_Friedgut} {Let $H=(V,E)$ be a hypergraph,
and $F_1,F_2,\ldots F_r$ subsets of $V$ such that every $v\in V$
belongs to at least $t$ of the sets $F_i$.  Let $H_i$ be the
projection hypergraphs: $H_i=(V,E_i),$ where $E_i=\{e\cap F_i : e\in
E\}$.  For each edge $e\in E$, define $e_i=e\cap F_i$, and assign
each $e_i$ a nonnegative real weight $w_i(e_i).$  Then
\[\Big(\sum_{e\in E}\BPi_{i=1}^{r}w_i(e_i)\Big)^t \leq \BPi_i\sum_{e_i\in E_i}w_i(e_i)^t\]
}
\end{theorem}

The first step in applying this theorem is to define appropriate
variables. Let $G=(V,E)$ be a $d$-regular graph, with its vertex set
$\{v_1,v_2, \ldots, v_N\}$. We will use $G$ to form an associated
matching hypergraph, $H=(E,{\cal M})$, where the vertex set of the
hypergraph is the edge set of $G$, and ${\cal M}$ is the sets of
matchings in $G$. Let $F_i$ be the set of edges incident to a vertex
$v_i\in V$.  Note that each edge in $E$ is covered twice by
$\bigcup_{i=1}^N F_i$, so we may take $t=2$.  We define the trace
sets, $E_i=\{F_i\cap m: m\in {\cal M}\},$ as the set of possible
intersections of a matching with the set of edges incident with
$v_i$. Let $m_i=m\cap F_i$. Then for all $i$, assign
\[w_i(m_i)= \left\{
\begin{array}{cc}
1&\mbox{ if } m_i=\emptyset \\
\sqrt{\lambda} &\mbox{ else}
\end{array}\right.\]
With these definitions we have $\sum_{m_i\in
E_i}w_i(m_i)^2=1+d\lambda$, and for a fixed $m$, $\BPi_i
w_i(m_i)=\sqrt{\lambda}^{(2|m|)}$. Putting these expressions into
Theorem \ref{weighted_entropy_Friedgut}, we have that
\[(\PFG)^2=\left(\sum_{m\in {\cal M}} \lambda^{|m|}\right)^2 \leq \BPi_{i=1}^{N}  (1+d\lambda).\]
Therefore,
$$
\label{PFG} \PFG \leq (1+d\lambda)^{\frac{N}{2}}.
$$
\qed

\begin{remark}
After the submission of this paper, L. Gurvits pointed out an alternative proof of Lemma \ref{pfg}, which applies to graphs with average degree $d$ and actually gives a slight improvement when $G$ does not have a perfect matching. By a result of Heilmann and Lieb \cite{HeilmannLieb}, the roots of $\PFG=0$ are all real and negative, and so we can write $\PFG = \prod_{i=1}^{\nu(G)} (1+\alpha_i \lambda)$ for some positive $\alpha_i$'s with $\sum \alpha_i = \left.(\PFG)'\right|_{\lambda=0}=|E(G)|=\frac{Nd}{2}$, where $\nu(G)$ is the size of the largest matching of $G$. Applying the arithmetic mean - geometric mean inequality to this expression we obtain
$$
\PFG \leq \left(1+\lambda\frac{\sum \alpha_i}{\nu(G)}\right)^{\nu(G)} = \left(1+\lambda\frac{Nd}{2\nu(G)}\right)^{\nu(G)} \leq \left(1+d\lambda\right)^{\frac{N}{2}}.
$$
\end{remark}
%
%

\noindent \textbf{Proof of Theorem \ref{dl}} We begin with the upper
bound. We may assume $0 < \ell < N/2$, since the extreme cases
$\ell=0, N/2$ are obvious. For fixed $\ell$, a single term of the
partition function $\PFG$ is bounded by the whole sum, and so by
Lemma \ref{pfg} we have $\ml(G) \lambda^{\ell} \leq \PFG \leq
(1+d\lambda)^\frac{N}{2}$ and
\begin{equation}
\label{PFGST} \ml(G) \leq (1+d\lambda)^\frac{N}{2}
\Big(\frac{1}{\lambda}\Big)^{\ell}.
\end{equation}
We take
$$
\lambda=\frac{\ell}{d\left(\frac{N}{2}-\ell\right)}
$$
to minimize the right hand
side of (\ref{PFGST}) and obtain the upper bound in Theorem \ref{dl}
(in the case $\ell=\frac{\alpha N}{2}$):
\begin{eqnarray*}
\log(\ml (G))
& \leq & \log \left(\frac{\frac{N}{2}}{\frac{N}{2}-\ell}\right)^{\frac{N}{2}}\left(\dfrac{d\big(\frac{N}{2}-\ell\big)}{\ell}\right)^{\ell}\\
& = & \frac{N}{2}\left(\frac{2\ell}{N}\log
d+H\left(2\ell/N\right)\right) \\
& = & \frac{N}{2}\left(\alpha\log d+H(\alpha)\right).
\end{eqnarray*}

We now turn to the lower bound. We begin by observing
\begin{equation} \label{matchings_in_DK}
\ml(\DK) = \sum_{a_1, \ldots a_{N/2d}:\atop{0 \leq a_i \leq d,~
\sum_i a_i = \ell}} \prod_{i=1}^{N/2d} {d \choose a_i}^2 a_i!
\end{equation}
Here the $a_i$'s are the sizes of the intersections of the matching
with each of the components of $\DK$, and the term ${d \choose
a_i}^2 a_i!$ counts the number of matchings of size $a_i$ in a
single copy of $K_{d,d}$. (The binomial term represents the choice
of $a_i$ endvertices for the matching from  each partition class,
and the factorial term tells us how many ways there are to pair the
endvertices from the top and bottom to form a matching.)

From Stirling's formula we have that there is an absolute constant
$c\geq 1$ such that for any $d \geq 1$ and $0 < a < d$,
\begin{equation} \label{stirling_bound}
\log \left({d \choose a}^2 a!\right) \geq a\log d + a \log
\frac{a}{d} -a\log e +2H(a/d)d -\log cd,
\end{equation}
and we may verify by hand that (\ref{stirling_bound}) holds also for
$a=0, d$. Combining (\ref{matchings_in_DK}) and
(\ref{stirling_bound}) we see that $\log (\ml(\DK))$ is bounded
below by
\begin{equation} \label{ais}
\frac{N}{2} \left(\frac{2\ell}{N}\log d - \frac{2\ell}{N}\log e
-\frac{\log cd}{d} + \frac{2}{N}\sum_{i=1}^{N/2d}\left(a_i\log
\frac{a_i}{d}+2H(a_i/d)d\right) \right)
\end{equation}
for any valid sequence of $a_i$'s. To get our lower bound in the
case $\ell = \alpha \frac{N}{2}$, we consider (\ref{ais}) for that
sequence of $a_i$'s in which each $a_i$ is either $\lfloor \alpha d
\rfloor$ or $\lceil \alpha d \rceil$. Note that by the mean value
theorem, there is a constant $c_\alpha
> 0$ such that both
$$
\log \frac{\lceil \alpha d \rceil}{d}, ~\log \frac{\lfloor \alpha d
\rfloor}{d} \geq \log \alpha - \frac{c_\alpha}{d}
$$
and
$$
H\left(\frac{\lceil \alpha d \rceil}{d}\right), ~H\left(
\frac{\lfloor \alpha d \rfloor}{d}\right) \geq H(\alpha) -
\frac{c_\alpha}{d}.
$$
(Here we use
$$
\left|\frac{\lceil \alpha d \rceil}{d} - \alpha\right|,~
\left|\frac{\lfloor \alpha d \rfloor}{d} - \alpha\right| \leq
\frac{1}{d}
$$
and $\alpha \neq 0, 1$.) Putting these bounds into (\ref{ais}) we
obtain
$$
\log (\ml(\DK)) \geq \frac{N}{2} \left(\alpha\log d + 2H(\alpha) +
\alpha\log \left(\frac{\alpha}{e}\right) +\Omega\left(\frac{\log
d}{d}\right)\right),
$$
with the constant in the $\Omega$ term depending on $\alpha$. \qed

\section{Counting Independent Sets}

In this section we prove the various assertions of Theorem
\ref{thm-indsets}. We begin with the second bound in
(\ref{ind_set_ub}). We use a result from \cite{Galvi2004}, which
states that for any $\gl > 0$ and any $d$-regular $N$-vertex
bipartite graph $G$, the weighted independent set  partition
function satisfies
\begin{equation} \label{fromGT}
Z_{\lambda}^{\rm{ind}}(G) :=\sum_{I \in \cI(G)} \gl^{|I|} \leq
\left(2(1+\gl)^d-1\right)^{\frac{N}{2d}}.
\end{equation}
Choose $\lambda$ so that $\frac{\lambda N}{2(1+\lambda)}=t$. Noting
that $i_t(G)\lambda^{\frac{\lambda N}{2(1+\lambda)}}$ is the
contribution to $Z_{\lambda}^{\rm{ind}}(G)$ from independent sets of
size $t$ we have
\begin{eqnarray}
\label{secondlast} i_t(G) & \leq &
\frac{Z_\lambda^{\rm{ind}}(G)}{\lambda^{\frac{\lambda
N}{2(1+\lambda)}}} \nonumber \\
& \leq &
\frac{\left(2(1+\gl)^d-1\right)^{\frac{N}{2d}}}{\lambda^{\frac{\lambda
N}{2(1+\lambda)}}} \label{usingGT} \\
& = & 2^\frac{N}{2d}
\left(\frac{1+\lambda}{\lambda^\frac{\lambda}{1+\lambda}}\right)^{N/2}\left(1-\frac{1}{2(1+\gl)^d-1}\right)^{\frac{N}{2d}}
\nonumber \\
& = & 2^{H\left(\frac{\lambda}{1+\lambda}\right)\frac{N}{2} +
\frac{N}{2d}}e^{-\frac{N}{4d(1+\gl)^d}}
\nonumber \\
& = & 2^{H\left(\frac{2t}{N}\right)\frac{N}{2} + \frac{N}{2d}-
\frac{N\log e }{4d}\left(1-\frac{2t}{N}\right)^d}. \nonumber
\end{eqnarray}
\noindent We use (\ref{fromGT}) to make the critical substitution in
(\ref{usingGT}).

To obtain the first bound in (\ref{ind_set_ub}) we need the
following analog of (\ref{fromGT}) for $G$ not necessarily
bipartite:
\begin{equation} \label{partition_function_non_bipartite}
Z_\lambda^{\rm{ind}}(G) \leq 2^{\frac{N}{d}}(1+\lambda)^\frac{N}{2}.
\end{equation}
From (\ref{partition_function_non_bipartite}) we easily obtain the
claimed bound, following the steps of the derivation of the second
bound in (\ref{ind_set_ub}) from (\ref{fromGT}). We prove
(\ref{partition_function_non_bipartite}) by using a more general
result on graph homomorphisms. For graphs $G=(V_1,E_1)$ and
$H=(V_2,E_2)$ set
$$
Hom(G,H)=\{f:V_1\rightarrow V_2~:~\{u,v\} \in E_1 \Rightarrow
\{f(u),f(v)\} \in E_2\}.
$$
That is, $Hom(G,H)$ is the set of graph homomorphisms from $G$ to
$H$. Fix a total order $\prec$ on $V(G)$. For each $v \in V(G)$,
write $P_\prec(v)$ for $\{w \in V(G):\{w,v\} \in E(G), w \prec v\}$
and $p_\prec(v)$ for $|P_\prec(v)|$. The following natural
generalization of a theorem of J. Kahn is due to D. Galvin (see
\cite{MadimanTetali} for a proof).
\begin{thm} \label{thm-main}
For any $d$-regular and $N$-vertex graph $G$ (not
necessarily bipartite) and any total order $\prec$ on $V(G)$,
$$
|Hom(G,H)| \leq \prod_{v \in V(G)}
|Hom(K_{p_\prec(v),p_\prec(v)},H)|^\frac{1}{d}.
$$
\end{thm}
If $G$ is bipartite with bipartition classes $\cE$ and $\cO$ and
$\prec$ satisfies $u \prec v$ for all $u \in \cE, v \in \cO$ then
Theorem \ref{thm-main} reduces to the main result of
\cite{Galvi2004}.

\medskip

To prove (\ref{partition_function_non_bipartite}), we first note
that (by continuity) it is enough to prove the result for $\lambda$
rational. Let $C$ be an integer such that $C\lambda$ is also an
integer, and let $H_C$ be the graph which consists of an independent
set of size $C\lambda$ and a complete looped graph on $C$ vertices,
with a complete bipartite graph joining the two. As described in
\cite{Galvi2004} we have, for any graph $G$ on $N$ vertices,
$$
|Hom(G,H_C)|=C^N Z_\lambda^{\rm{ind}}(G).
$$
For $G$ $d$-regular and $N$-vertex, we apply Theorem \ref{thm-main}
twice to obtain
\begin{eqnarray*}
Z_\lambda^{\rm{ind}}(G) & = & \frac{|Hom(G,H_C)|}{C^N} \\
& \leq & \frac{\prod_{v \in V(G)}
|Hom(K_{p_\prec(v),p_\prec(v)},H_C)|^\frac{1}{d}}{C^N} \\
& = & \frac{\prod_{v \in V(G)}
\left(C^{2p_\prec(v)}Z_\lambda^{\rm{ind}}(K_{p_\prec(v),p_\prec(v)})\right)^\frac{1}{d}}{C^N}
\\
& \leq & \frac{C^{\frac{2\sum_{v \in V(G)} p_\prec(v)}{d}} \prod_{v
\in
V(G)} \left(2(1+\lambda)^{p_\prec(v)}\right)^\frac{1}{d}}{C^N} \\
& = & 2^\frac{N}{d} \frac{C^{\frac{2\sum_{v \in V(G)}
p_\prec(v)}{d}} (1+\lambda)^{\frac{\sum_{v \in V(G)}
p_\prec(v)}{d}}}{C^N}.
\end{eqnarray*}
Now noting that
$$
\sum_{v \in V(G)} p_\prec(v) = |E(G)| = \frac{Nd}{2}
$$
we obtain
$$
Z_\lambda(G) \leq 2^{\frac{N}{d}}(1+\lambda)^\frac{N}{2},
$$
as claimed.

We now turn to the third bound in (\ref{ind_set_ub}). Fix a perfect
matching of $G$ joining a set of vertices $A \subseteq V(G)$ of size
$N/2$ to the set $B:=V(G)\setminus A$. Let $f$ be the bijection from
subsets of $A$ to subsets of $B$ that moves the set along the chosen
matching. Every independent set in $G$ of size $t$ is of the form
$I_A \cup I_B$ where $I_A \subseteq A$, $I_B \subseteq B$, $f(A)
\cap B = \emptyset$ and $|A|+|B|=t$. We therefore count all the
independent sets of size $t$ (and more) by choosing a subset of $A$
of size $t$ (${N/2 \choose t}$ choices) and a subset of this set to
send to $B$ via $f$ ($2^t$ choices).

To obtain the first bound in (\ref{ind_set_lb}), we introduce a
probabilistic framework and use Markov's inequality.  If we divide a
set of size $N/2$ into $N/2d$ blocks of size $d$ and choose a
uniform subset of size $t$, then the probability that this set
misses a particular block is ${N/2-d \choose t}/{N/2 \choose t}$.
Let $X$ be a random variable representing the number of blocks that
the $t$-set misses.  Let $b_k$ equal the number of $t$-sets which
miss exactly $k$ blocks. Then $\mathbb{P}(X=k)=b_k/\binom{N/2}{t}$.
Let $\chi_A$ be the indicator variable for the event $A$.  Then
$$
X = \sum_{i=0}^{\frac{N}{2d}}\chi_{\{\rm{block}~i~\rm{empty}\}}
$$
and by linearity of expectation the expected number of blocks missed
satisfies
\begin{equation} \label{bound_on_mu}
\mu := \mathbb{E}(X)=\frac{N}{2d}\frac{{\frac{N}{2}-d \choose
t}}{{\frac{N}{2} \choose t}} \leq
\frac{N}{2d}\left(1-\frac{2t}{N}\right)^d.
\end{equation}

From Markov's inequality we have
$$
\sum_{k=0}^{c\mu}\mathbb{P}(X=k)=\mathbb{P}(X\leq c\mu) \geq
\left(1-\frac{1}{c}\right).
$$
We substitute the previously discussed value for $\mathbb{P}(X=k)$,
yielding the inequality
\begin{equation} \label{markov}
\sum_{k=0}^{c\mu}b_k\geq
\left(1-\frac{1}{c}\right)\binom{\frac{N}{2}}{t}.
\end{equation}

How many independent sets of size $t$ does $\DK$ have? To choose an
independent set from $\DK$ of size $t$, we first create a
bipartition $\cE \cup \cO$ of $\DK$ by choosing (arbitrarily) one of
the bipartition classes of each of the $N/2d$ $K_{d,d}$'s of $\DK$
to be in $\cE$. We then choose a subset of $\cE$ of size $t$. The
number of subsets of $\cE$ which have empty intersection with
exactly $k$ of the $K_{d,d}$'s that make up $\DK$ is precisely
$b_k$. Each of these subsets corresponds to $2^{\frac{N}{2d} - k}$
independent sets in $\DK$. Combining this observation with
(\ref{bound_on_mu}) and (\ref{markov}) we obtain the first bound in
(\ref{ind_set_lb}):
\begin{eqnarray*}
i_t(\DK)  & =  & 2^\frac{N}{2d} \sum_{k\geq 0} 2^{-k}b_k \\
& \geq & 2^{\frac{N}{2d}-c\mu} \sum_{k = 0}^{c\mu} b_k \\
& \geq &
\left(1-\frac{1}{c}\right)\binom{\frac{N}{2}}{t}
2^{\frac{N}{2}\left(\frac{1}{d}-\frac{c}{d}\left(1-\frac{t}{M}\right)^d\right)}.
\end{eqnarray*}

Finally we turn to the second bound in (\ref{ind_set_lb}). We obtain
the claimed bound by considering all of the independent sets
whose intersection with each component of $\DK$ has size either $0$ or
$1$:
$$
i_t(\DK) \geq (2d)^t {\frac{N}{2d} \choose t}.
$$
After a little algebra, the right hand side above is seen to be
exactly the right hand side of the second bound in
(\ref{ind_set_lb}).

\end{document}